\documentclass[12pt]{article}
\usepackage{amssymb}
\newtheorem{theorem}{Theorem}

\newtheorem{proposition}[theorem]{Proposition}
\newtheorem{lemma}[theorem]{Lemma}
\newtheorem{corollary}[theorem]{Corollary}
\newtheorem{remark}[theorem]{Remark}
\newtheorem{remarks}[theorem]{Remarks}

\newcommand{\R}{\mathbb{R}}

\newcommand{\hess}{\mbox{Hess\,}}

\newcommand{\E}{{\cal E}}
\newcommand{\po}{{\hspace*{-1ex}}{\bf .  }}

\newcommand{\su}{submanifold }

\def\<{\langle}

\def\>{\rangle}

\def\va{\varphi}

\def\bea{\begin{eqnarray*} }
\def\eea{\end{eqnarray*} }
\def\be{\begin{equation} }
\def\ee{\end{equation} }

\def\proof{\noindent{\it Proof: }}
\def\qed{\ifhmode\unskip\nobreak\fi\ifmmode\ifinner
\else\hskip5 pt \fi\fi\hbox{\hskip5 pt \vrule width4 pt
height6 pt  depth1.5 pt \hskip 1pt }}
\begin{document}
\title{\bigskip
\bigskip
Conformal de Rham decomposition of Riemannian manifolds.}
\author{Ruy Tojeiro}
\date{}
\maketitle
\noindent {\bf Abstract:} {\small We prove conformal versions of
the local decomposition theorems of de Rham and Hiepko of a
Riemannian manifold as a  Riemannian or a warped product of
Riemannian manifolds.   Namely, we give necessary and sufficient
conditions for a Riemannian manifold to be locally conformal to
either a Riemannian or a warped product. We also obtain other
related de Rham-type decomposition theorems. As an application, we
study Riemannian manifolds that admit a Codazzi tensor with two
distinct
 eigenvalues everywhere. }
\vspace{3ex}

\vspace{3ex}

\vspace{3ex}


\bigskip
\noindent {\bf\large  \S 1. Introduction. }\vspace{3ex}

An important initial step in understanding a complicated
mathematical object is to decompose it into simpler "irreducible"
components. In differential geometry, a fundamental result in this
direction is the decomposition theorem of de Rham, which gives
necessary and sufficient conditions for a Riemannian manifold to
split, both locally and globally, into a Riemannian product of
Riemannian manifolds \cite{Rh}. Several  significant
generalizations of de Rham's theorem have been  obtained. For
instance, a similar  characterization was given by Hiepko
\cite{hi} of Riemannian manifolds that split as a warped product
of Riemannian manifolds. More recently, a suitable setting for
treating such decomposition results was introduced in \cite{rs} by
defining the notions of {\it netted manifolds\/} and of {\it net
morphisms\/} between them.

      A {\it net\/} on a connected $C^{\infty}$-manifold $M$ is a splitting
      $TM=\oplus_{i=1}^k E_i$ of the tangent bundle of $M$ by a family of  integrable subbundles. If $M$ is a Riemannian manifold and the subbundles are mutually orthogonal then the net is said to be an {\it orthogonal net\/}. The canonical net on a product manifold
$M=\Pi_{i=1}^k M_i$ is called the {\it product net\/}. A
$C^\infty$-map $\psi\colon\;M\to N$ between two {\it netted
manifolds\/}  $(M,{\cal E})$, $(N,{\cal F})$, that is,
$C^\infty$-manifolds $M,N$ equipped with nets ${\cal
E}=(E_i)_{i=1,\ldots,k}$ and ${\cal F}=(F _i)_{i=1,\ldots,k}$,
 respectively, is called a {\it net morphism\/} if $\psi_*E_i(p)\subset F_i(\psi(p))$
 for all $p\in M$, $1\leq i\leq k$,
or equivalently, if for any $p\in M$ the restriction
$\psi|_{L_i^{{\cal E}}(p)}$ to the leaf of $E_i$ through $p$ is a
$C^\infty$-map into the leaf $L_i^{{\cal F}}(\psi(p))$ of $F_i$
through $\psi(p)$. The net morphism $\psi$ is said to be a {\it
net isomorphism\/} if, in addition, it is a diffeomorphism and
$\psi^{-1}$ is also a net morphism. A net ${\cal E}$ on $M$ is
said to be {\em locally decomposable\/} if for every point $p\in
M$ there exist a neighborhood $U$ of $p$ in $M$ and a net
isomorphism $\psi$ from $(U,{\cal E}|_U)$ onto a product manifold
$\Pi_{i=1}^k M_i$. The map $\psi^{-1}\colon\;\Pi_{i=1}^k M_i\to U$
is called a {\em product representation\/} of $(U,{\cal E}|_U)$.

     A general problem in this context is to determine necessary and sufficient conditions
     for an orthogonal net to admit  locally a product representation whose induced metric
     is of some nice particular type, for instance a Riemannian product or a warped product
     of Riemannian metrics on the factors as in de Rham and Hiepko theorems, respectively.

    In this article we consider this problem from a conformal point of view.
    Our main results are conformal versions of the local decomposition theorems
    of de Rham and Hiepko. Namely, we find  necessary and sufficient conditions
    for an orthogonal net to admit locally a product representation whose induced
    metric is {\it conformal\/} to either  a Riemannian product or a warped product
    of Riemannian metrics on the factors. We also solve the above
    problem for other natural types of metrics and orthogonal
    nets.

     As an application, we consider the problem of determining the restrictions that
     are imposed on a Riemannian manifold by the existence of a Codazzi tensor with
exactly two distinct
      eigenvalues (cf.   \cite{de} and \cite{ds}).
      We start from
      the observation that the orthogonal net  determined by the eigenbundles of such
      a  Codazzi tensor admits locally a product representation whose induced metric
      is a twisted product of  Riemannian metrics on the factors, and study further
      properties of the twisted product metric under additional assumptions on the
      Codazzi tensor. In particular, we determine all Riemannian
      manifolds that carry Codazzi tensors with  exactly two distinct
      eigenvalues, both of which are constant along the
      corresponding eigenbundles, as well as the tensors
      themselves.

  In a forthcoming paper \cite{to} we study the problem posed by Burstall \cite{bu}
  of developing a theory of isothermic Euclidean submanifolds of higher dimensions
  and codimensions. Recall that a surface  in Euclidean three-space is called
  {\it isothermic\/} if, away from umbilic points,  its curvature lines form an
  isothermic net, that is, there exist locally conformal coordinates that diagonalize
  the second fundamental form. Studying conformal decomposition theorems of Riemannian
  manifolds was in part motivated by the problem of looking for a suitable extension
  of the intrinsic notion of an isothermic net. In fact, our conformal version
  of the local
  de Rham theorem can be seen as a far reaching generalization of a classical
  characterization of isothermic nets of curves in terms of their geodesic curvatures \cite{da}.

     To conclude this introduction, we point out that the results of this paper,
     as well as their proofs, remain valid for pseudo--Riemannian manifolds.\vspace{2ex}

\noindent {\bf\large  \S 2. Twisted and warped products.
}\vspace{2ex}

    In this section we recall from \cite{mrs} some basic definitions and results on warped
    and twisted products.
If $M=\Pi_{i=0}^k M_i$ is the product of $C^\infty$-manifolds
$M_0,\ldots,M_k$, then $\<\;,\;\>$ is called a {\em twisted
product metric\/} on $M$ if there exist Riemannian metrics
$\<\;,\;\>_i$ on $M_i$, $0\leq i\leq k$, and a $C^\infty$ {\em
twist-function \/} $\rho=(\rho_0,\ldots,\rho_k)\colon\;M\to
\R^{k+1}_+$ such that
$$\<\;,\;\>=\sum_{i=0}^k\rho_i^2\pi_i^*\<\;,\;\>_i,$$
where $\pi_i\colon\,M\to M_i$ denotes the canonical projection.
Then $(M,\<\;,\;\>)$ is said to be a {\em twisted product\/}  and
is denoted by $^\rho\Pi_{i=0}^k (M_i,\<\;,\;\>_i)$. When
$\rho_1,\ldots,\rho_k$ are independent of $M_1,\ldots, M_k$, that
is, there exist $\tilde{\rho}_i\in C^\infty(M_0)$ such that
$\rho_i=\tilde{\rho}_i\circ \pi_0$ for $i=1,\ldots, k$ and, in
addition, $\rho_0$ is identically one,  then $\<\;,\;\>$ is called
a {\em warped product metric\/} and
$(M,\<\;,\;\>):=(M_0,\<\;,\;\>_0)\times_{\rho}\Pi_{i=1}^k
(M_i,\<\;,\;\>_i)$ a {\em warped product\/} with {\em warping
function\/} $\rho=(\rho_1,\ldots,\rho_k)$. If $\rho_i$ is
identically $1$ for  $i=0,\ldots, k$, the metric $\<\;,\;\>$ is
the usual {\em Riemannian product metric\/}, in which case
$(M,\<\;,\;\>)$ is called a {\em Riemannian product\/}.

 The next result from \cite{mrs} relates the Levi-Civita connections of a
 twisted product metric and the corresponding  Riemannian product metric on a
  product manifold.

\begin{proposition} $\mbox{{\bf \cite{mrs}}}$ \label{prop:cons} Let
$(M,\<\;,\;\>)=^\rho \Pi_{i=0}^k (M_i,\<\;,\;\>_i)$ be a twisted product with
twist function $\rho=(\rho_0,\ldots,\rho_k)$ and product net ${\cal E}=(E_i)_{i=0,\ldots,k}$,
let $\nabla$ and $\tilde{\nabla}$ be the Levi-Civita connections of $\<\;,\;\>$
and of the product metric $\<\;,\;\>^\sim$, respectively,
and let $U_i=-\nabla (\log\circ \rho_i)$, $0\leq i\leq k$,
where the gradient is calculated with respect to $\<\;,\;\>$. Then
\be \label{eq:nablas}
\nabla_X Y=\tilde{\nabla}_X Y +\sum_{i=0}^k (\<X^i,Y^i\>U_i - \<X,U_i\>Y^i- \<Y,U_i\>X^i).
\ee
\end{proposition}

    Here and throughout the paper, if ${\cal E}=(E_i)_{i=0,\ldots,k}$
    is an orthogonal net on a Riemannian manifold then writing a vector field with a
    superscript $i$ (resp., $\perp_i$) indicates taking the $E_i$-component
    (resp., $E_i^\perp$-component) of that vector field. Moreover, we always denote
    sections of $E_i$ (resp., $E_i^\perp$) by $X_i$ and $Y_i$ (resp., $X_{\perp_i}$
    and $Y_{\perp_i}$).

An orthogonal net ${\cal E}=(E_i)_{i=0,\ldots,k}$ on a Riemannian manifold  $M$ is
called a $TP$-{\em net\/} if $E_i$ is umbilical and $E_i^\perp$ is integrable for
 $i=0,\ldots, k$.
Recall that a subbundle $E$ of $TM$ is {\em umbilical\/} if there exists a vector
field $\eta$ in $E^\perp$ such that
$$\<\nabla_X Y,Z\>=\<X,Y\>\<\eta,Z\>\;\;\mbox{for all }\;X,Y\in E, \;Z\in E^\perp.$$
The vector field $\eta$ is called the {\em mean curvature normal\/} of $E$. If, in addition,
$$\<\nabla_X \eta,Z\>=0\;\;\mbox{for all }\;X\in E, \;Z\in E^\perp,$$
then $E$ is said to be {\em spherical\/}. If $E$ is umbilical and its mean curvature normal
vanishes identically, then it is called {\it totally geodesic\/} (or {\it auto-parallel\/}).
An umbilical distribution is automatically integrable, and the leaves are umbilical
submanifolds of $M$. When $E$ is totally geodesic or spherical, then its leaves are
totally geodesic or spherical submanifolds, respectively. By a {\it spherical submanifold\/}
we mean an umbilical \su whose mean curvature vector is parallel with respect to the
normal connection.

    An orthogonal net ${\cal E}=(E_i)_{i=0,\ldots,k}$ is called a $WP$-{\em net\/}
    if $E_i$ is spherical and $E_i^\perp$ is totally geodesic  for $i=1,\ldots, k$ .
    In a $WP$-net the subbundle $E_0$ is automatically totally geodesic and $E_0^\perp$
    is integrable (see Proposition $3$ in \cite{mrs} or Lemma~\ref{prop:cwp} below);
    in particular, every $WP$-net is also a $TP$-net.

     It follows from  (\ref{eq:nablas}) that
\be\label{eq:1a}
(\nabla_{X_i}Y_i)^{\perp_i}=\<X_i,Y_i\>U_i^{\perp_i}, \ee and
\be\label{eq:1b} (\nabla_{X_{\perp_i}}Y_{\perp_i})^i=\sum_{j\neq
i}\<X_{\perp_i}^j,Y_{\perp_i}^j\>U_j^i. \ee Equation (\ref{eq:1a})
implies that $E_i$ is umbilical  with mean curvature normal
$H_i=U_i^{\perp_i}$. Thus the  product net of a twisted product is
a $TP$-net. For a warped product, we have that $H_i=U_i$ for
$i=1,\ldots, k$, because $\rho_i$ depends only on $M_0$,  and that
$E_i^\perp$ is totally geodesic, as follows from (\ref{eq:1b}).
Then $\<\nabla_{X_i}H_i,X_{\perp_i}\>=
\<\nabla_{X_{\perp_i}}H_i,X_i\>=0,$ where the first equality holds
because $H_i$ is a gradient vector field. Thus $E_i$ is spherical,
and hence ${\cal E}$ is a $WP$-net. The converses also hold (see
Proposition $4$ of \cite{mrs}):

\begin{proposition}\label{prop:tpwp} $\mbox{{\bf \cite{mrs}}}$
On a connected product manifold $M=:\Pi_{i=0}^k M_i$ the product
net ${\cal E}=(E_i)_{i=0,\ldots,k}$ is a $TP$-net (resp.,
$WP$-net) with respect to a Riemannian metric $\<\;,\;\>$ on $M$
if and only if $\<\;,\;\>$  is a twisted (resp., warped) product
 metric on $M$.
\end{proposition}

The following result from \cite{mrs} (see Corollary $1$ of \cite{mrs})
contains the local version of Hiepko's  decomposition theorem \cite{hi}.

\begin{theorem} \label{cor:isodrh} $\mbox{{\bf \cite{mrs}}}$
Let ${\cal E}=(E_i)_{i=0,\ldots,k}$ be a $TP$-net (resp., $WP$-net)
on a Riemannian manifold $M$. Then, for every point $p\in M$ there exists a
local product representation $\psi\colon\;\Pi_{i=0}^k M_i\to U$ of ${\cal E}$
with $p\in U\subset M$, which is an isometry with respect to a twisted  product
(resp., warped product) metric on $\Pi_{i=0}^k M_i$.
\end{theorem}

    Theorem \ref{cor:isodrh} is a consequence of  Proposition \ref{prop:tpwp}
    and the following basic criterion for local decomposability of a net on an
    arbitrary $C^\infty$-manifold (cf. Theorem $1$ of \cite{rs}).
\begin{proposition} \label{prop:dec} $\mbox{{\bf \cite{rs}}}$
A net ${\cal E}=(E_i)_{i=0,\ldots,k}$ on a $C^\infty$-manifold is
locally decomposable if and only if $E_i^{\perp}:=\oplus_{j\neq i}
E_i$ is integrable for $i=0,\ldots, k$.
\end{proposition}
\vspace{3ex}

\noindent {\bf\large  \S 3. Quasi-warped products.}\vspace{2ex}

  We say that a Riemannian metric  on a product manifold $M=\Pi_{i=0}^k M_i$  is a
  {\it quasi-warped product metric\/} if it is a twisted product metric  with twist
  function $\rho=(\rho_0,\ldots,\rho_k)$ and, in addition, $\rho_0$ is identically $1$
  and $\rho_i$  depends only on $M_0$ and $M_i$ for  $i=1,\ldots, k$.
  In this section we characterize the orthogonal nets that admit locally a
  product representation whose induced metric is either isometric or conformal
  to a quasi-warped product metric. We start with  a few preliminary facts.

\begin{lemma}\label{prop:cwp} Let ${\cal E}=(E_i)_{i=0,\ldots,k}$ be an orthogonal
net on a Riemannian manifold such that $E_i$ and $E_i^\perp$ are
umbilical for  $i=1,\ldots, k$. Then  $E_0^\perp$ is integrable
and $E_0$ is umbilical with mean curvature normal
$H_0=\sum_{i=1}^k \eta_i$, where $\eta_i$ is the mean curvature
normal of $E_i^\perp$. Therefore ${\cal E}$ is a $TP$-net.
Moreover, if $E_i^\perp$ is spherical for $i=1,\ldots, k$ then the
same holds for $E_0$.
\end{lemma}
\proof  First notice that $[X_i,Y_i]\in E_i\subset E_0^\perp$ for
$i=1,\ldots, k$, because $E_i$ is umbilical, and hence integrable.
Now, using that $E_i^\perp$ is umbilical with mean curvature
normal $\eta_i$ for $i=1,\ldots, k$ we have for all $i,j=1,\ldots,
k$ with $i\neq j$  that \be\label{eq:ineqj} \<\nabla_{X_j} X_i,
X_0\>=-\<X_i,\nabla_{X_j} X_0\>=-\<X_j,X_0\>\<\eta_i,X_i\>=0,\ee
thus $\nabla_{X_j} X_i\in E_0^\perp$. Therefore $E_0^\perp$ is
integrable. That $E_0$ is umbilical with mean curvature normal
$H_0=\sum_{i=1}^k \eta_i$ follows from
$$(\nabla_{X_0} Y_0)^{\perp_0}=\sum_{i=1}^k (\nabla_{X_0} Y_0)^{i}=
\sum_{i=1}^k \<X_0,Y_0\>\eta_i=\<X_0,Y_0\>\sum_{i=1}^k \eta_i.$$
Finally, since $E_j^\perp$ is umbilical we have that
$(\nabla_{X_0}\eta_i)^j=0$ for $i,j=1,\ldots,k$ with $i\neq j$,
and hence
$$ (\nabla_{X_0} H_0)^{\perp_0}=\sum_{i=1}^k (\nabla_{X_0} \eta_i)^{\perp_0}=
\sum_{i=1}^k (\nabla_{X_0} \eta_i)^i.$$
We conclude that $E_0$ is spherical whenever $E_i^\perp$ is spherical for  $i=1,\ldots, k$.\qed

\begin{lemma} \label{le:rhos} Let $(M,\<\;,\;\>)=\;^\rho\Pi_{i=0}^k (M_i,\<\;,\;\>_i)$
be a twisted product  and let ${\cal E}=(E_i)_{i=0,\ldots,k}$ be
its product net. Then for every fixed $i\in\{0,\ldots,k\}$ the
following are equivalent:
\begin{itemize}
\item[{\em $(i)$}]  $\rho_j/\rho_\ell$ does not depend on $M_i$ for all
$j,\ell\in \{0,\ldots,k\}$ with $j,\ell\neq i$;
\item[{\em $(ii)$}] $H_j^i=H_\ell^i$ for all $j,\ell\in \{0,\ldots,k\}$ with $j,\ell\neq i$,
where $H_j=U_j^{\perp_j}$ ($U_j=-\nabla \log\circ \rho_j)$, is the
mean curvature normal of $E_j$;
 \item[{\em $(iii)$}]  $E_i^\perp$ is umbilical.
\end{itemize}
If these assertions are true, then the mean curvature normal
$\eta_i$ of $E_i^\perp$ coincides with the vectors $H_j^i$, $j\neq
i$.
\end{lemma}
\proof  We have that $H_j^i-H_\ell^i=-(\nabla (\log\circ
(\rho_j/\rho_\ell)))^i$, which gives the equivalence between $(i)$
and $(ii)$. The equivalence between $(ii)$ and $(iii)$ follows
from (\ref{eq:1b}). \qed

\begin{proposition} \label{prop:cqwp} Let $M=\Pi_{i=0}^k M_i$ be a connected product
manifold and   ${\cal E}=(E_i)_{i=0,\ldots,k}$  its product net. Then $E_i$ is umbilical
and $E_i^\perp$ is totally geodesic (resp., $E_i$ and $E_i^\perp$ are umbilical)
for  $i=1,\ldots, k$ with respect to a Riemannian metric $\<\;,\;\>$ on $M$
if and only if $\<\;,\;\>$ is (resp., $\<\;,\;\>$ is conformal to) a
quasi-warped product metric.
\end{proposition}
\proof Assume that  $\<\;,\;\>=\va^2\<\;,\;\>^\sim$ is conformal
to a quasi-warped product metric
$\<\;,\;\>^\sim=\pi_0^*\<\;,\;\>_0+\sum_{i=1}^k\tilde{\rho}_i^2\pi_i^*\<\;,\;\>_i$.
Then it is a twisted product metric with twist function
$\rho=(\va,\va\tilde{\rho}_1,\ldots,\va\tilde{\rho}_k)$. Therefore
$E_i$ is umbilical for $i=1,\ldots, k$ (in fact for $i=0,\ldots,
k$) and, by Lemma \ref{le:rhos}, the same holds for $E_i^\perp$.
On the other hand, if $\<\;,\;\>$ is a quasi-warped product
metric, then (\ref{eq:1b}) implies that $E_i^\perp$ is totally
geodesic for $i=1,\ldots,k$.

        Conversely, assume that $E_i$ and $E_i^\perp$ are umbilical
        for  $i=1,\ldots, k$ with respect to a Riemannian metric $\<\;,\;\>$ on $M$.
        Then ${\cal E}$ is a $TP$-net by  Lemma \ref{prop:cwp}, and hence $\<\;,\;\>$ is a
        twisted product metric
$\<\;,\;\>=\sum_{i=0}^k\rho_i^2\pi_i^*\<\;,\;\>_i$ by Proposition
\ref{prop:tpwp}. Moreover, by Lemma \ref{le:rhos} we have that
$\tilde{\rho}_i=\rho_i/\rho_0$ only depends on $M_i$ and $M_0$ for
$i=1,\ldots, k$. Thus
$\<\;,\;\>^\sim=\pi_0^*\<\;,\;\>_0+\sum_{i=1}^k\tilde{\rho}_i^2\pi_i^*\<\;,\;\>_i$
 is a quasi-warped product metric and  $\<\;,\;\>=\rho_0^2\<\;,\;\>^\sim$.
 If, in addition, $E_i^\perp$ is totally geodesic for $i=1,\ldots, k$,
 then (\ref{eq:1b}) implies that $\rho_0$ only depends on $M_0$,
 and hence $\<\;,\;\>=\pi_0^*\<\;,\;\>^*_0+\sum_{i=1}^k\rho_i^2\pi_i^*\<\;,\;\>_i$
 is a quasi-warped product metric, where $\<\;,\;\>^*_0=\tilde{\rho}_0^{2}\<\;,\;\>_0$
 for $\tilde{\rho}_0\in C^{\infty}(M_0)$ such that $\tilde{\rho}_0\circ \pi_0=\rho_0$.\qed

\begin{remark} \label{re:unet}{\em In Proposition \ref{prop:cqwp} we have  that $E_i$
is umbilical for $i=1,\ldots, k$ and that $E_i^\perp$ is umbilical
for $i=0,\ldots, k$ (and hence also $E_0$ is umbilical by Lemma
\ref{prop:cwp}) with respect to a Riemannian metric $\<\;,\;\>$ on
$M$ if and only if $\<\;,\;\>$ is conformal to a quasi-warped
product metric
$\<\;,\;\>^\sim=\pi_0^*\<\;,\;\>_0+\sum_{i=1}^k\tilde{\rho}_i^2\pi_i^*\<\;,\;\>_i$
and, in addition, $\tilde{\rho}_i/\tilde{\rho}_j$ does not depend
on $M_0$ for $i,j=1,\ldots, k$,
 or equivalently,  there exists $Y_0\in E_0$ (not necessarily a gradient
 vector field) such that
 $(\nabla \log\circ \tilde{\rho}_i)^0=Y_0$ for $i=1,\ldots, k$. This follows by applying
 Lemma \ref{le:rhos} once more for $i=0$.}\end{remark}

\begin{theorem}\po \label{prop:qwp} Let ${\cal E}=(E_i)_{i=0,\ldots,k}$
be an orthogonal net on a Riemannian manifold~$M$. Assume that
$E_i$ is umbilical and $E_i^\perp$ is totally geodesic (resp.,
$E_i$ and $E_i^\perp$ are umbilical) for  $i=1,\ldots, k$.
Then for every point $p\in M$ there exists a local product
representation $\psi\colon\;\Pi_{i=0}^k M_i\to U$ of ${\cal E}$
with $p\in U\subset M$, which is an isometry (resp., a conformal
diffeomorphism) with respect to a quasi-warped product metric on
$\Pi_{i=0}^k M_i$.
\end{theorem}
\proof By Lemma  \ref{prop:cwp}, the net ${\cal E}$ is a $TP$-net, and thus it is locally
decomposable by Proposition \ref{prop:dec}. Let
$\psi\colon\;\Pi_{i=0}^k M_i\to U$ be a local product representation
with $p\in U\subset M$
and let  $\Pi_{i=0}^k M_i$ be endowed with the metric induced by $\psi$.
Then the product net of $\Pi_{i=0}^k M_i$ satisfies that same properties as ${\cal E}$
and Proposition \ref{prop:cqwp}  concludes the proof.\vspace{3ex}\qed

\noindent {\bf\large  \S 4. Local conformal versions of de Rham and Hiepko theorems. }\vspace{2ex}

   In this section we prove our main results. Namely, we derive conformal versions of the
   local  decomposition theorems of de Rham and Hiepko and prove some related decomposition
   results.  \vspace{.5ex}

  We say that an orthogonal net  ${\cal E}=(E_i)_{i=0,\ldots,k}$ on a Riemannian manifold
  is a  {\em conformal warped product net\/},  or a $CWP$-net for short, if for $i=1,\ldots, k$
  it holds that
\be\label{eq:hs} E_i \mbox{ and } E_i^\perp \mbox{ are umbilical and }\<\nabla_{X_{\perp_i}}\eta_i,X_i\>=\<\nabla_{X_i}H_i,X_{\perp_i}\>,\ee
where $H_i$ and $\eta_i$ are the mean curvature normals of $E_i$ and $E_i^\perp$, respectively.
 If, in addition, also $E_0^\perp$ is umbilical, then we say that ${\cal E}$ is a
 {\em conformal product net\/}, or a $CP$-net for short. We first observe a few
  elementary facts on $CP$-nets and $CWP$-nets.

\begin{proposition} \label{prop:sph}
\begin{itemize} \item[{\em $(i)$}] For a $CP$-net ${\cal E}=(E_i)_{i=0,\ldots,k}$ condition
(\ref{eq:hs})  holds also for $i=0$;
\item[{\em $(ii)$}] If ${\cal E}=(E_i)_{i=0,\ldots,k}$ is an orthogonal net such that $E_i$ and $E_i^\perp$ are spherical for  $i=1,\ldots, k$ then it is a  $CWP$-net.  If, in addition, also $E_0^\perp$ is spherical then it is a  $CP$-net.
\item[{\em $(iii)$}] If ${\cal E}=(E_i)_{i=0,\ldots,k}$ is a $CP$-net (resp., a $CWP$-net), then for each $i=0,\ldots, k$ (resp., $i=1,\ldots, k$ ) one of $E_i$ and $E_i^\perp$ being spherical implies the same for the other.
\end{itemize}
\end{proposition}
\proof If ${\cal E}=(E_i)_{i=0,\ldots,k}$ is a $CP$-net, then we
have from Lemma \ref{prop:cwp} that $E_0$ is umbilical with mean
curvature normal $H_0=\sum_{i=1}^k \eta_i$, where $\eta_i$ is the
mean curvature normal of $E_i^\perp$. Thus, in order to prove
$(i)$ it remains to verify that
$\<\nabla_{X_{\perp_0}}\eta_0,X_0\>=\<\nabla_{X_0}H_0,X_{\perp_0}\>$,
where $\eta_0$ is the mean curvature normal of $E_0^\perp$. First
recall that $H_i^0=\eta_0$ for $i=1,\ldots,k$ by Lemma
\ref{le:rhos}, where $H_i$ is the mean curvature normal of $E_i$.
 Then,
$$\begin{array}{l}\<\nabla_{X_0}H_0,X_{\perp_0}\>=\sum_{i=1}^k \<\nabla_{X_0}\eta_i,X_{\perp_0}\>=\sum_{i=1}^k \<\nabla_{X_0}\eta_i,X_{\perp_0}^i\>=\sum_{i=1}^k \<\nabla_{X_{\perp_0}^i}H_i,X_0\>\\
\hspace{14.6ex}=\sum_{i=1}^k
\<\nabla_{X_{\perp_0}^i}H_i^0,X_0\>=\sum_{i=1}^k
\<\nabla_{X_{\perp_0}^i}\eta_0,X_0\>=\<\nabla_{X_{\perp_0}}\eta_0,X_0\>,\end{array}$$
where in the second equality we have used that
$(\nabla_{X_0}\eta_i)^j=0$ for $i,j=1,\ldots,k$ with $i\neq j$
because $E_j^\perp$ is umbilical and in the fourth one that
$\nabla_{X_{\perp_0}^i}H_i^j=0$ for $i,j=1,\ldots,k$ with $i\neq
j$ because $E_0^\perp$ is umbilical.  Clearly, condition
(\ref{eq:hs}) is satisfied if both $E_i$ and $E_i^\perp$ are
spherical. Moreover, if (\ref{eq:hs}) holds  then one of $E_i$ or
$E_i^\perp$ being spherical implies the same for the other. Taking
into account $(i)$ and the last statement in Lemma \ref{prop:cwp},
all the assertions in $(ii)$ and $(iii)$ follow.\vspace{1ex}\qed

\begin{proposition} \label{prop:conf} On a connected and simply connected
product manifold $M=\Pi_{i=0}^k M_i$
the product net ${\cal E}=(E_i)_{i=0,\ldots,k}$ is a $CWP$-net
(resp., $CP$-net) with respect to a Riemannian metric $\<\;,\;\>$
on $M$ if and only if $\<\;,\;\>$ is conformal to a warped product
metric (resp., to a Riemannian product metric) on $M$.
\end{proposition}
\proof If $\<\;,\;\>^\sim=\pi_0^*\<\;,\;\>^\sim_0+\sum_{i=1}^k
\tilde{\rho}_i^2\pi_i^*\<\;,\;\>^\sim_i$ is a warped product
metric on $M$ and  $\<\;,\;\>=\va^2\<\;,\;\>^\sim$ for some
$\va\in C^\infty(M)$ then $\<\;,\;\>$ is a twisted product metric
with twist function
$\rho=(\rho_0,\ldots,\rho_k)=(\va,\va\tilde{\rho}_1,\ldots,\va\tilde{\rho}_k)$.
Set $W=-\nabla \log\circ\va$, $\tilde{U}_i =-\nabla \log\circ
\tilde{\rho}_i\in E_0$ and $U_i=-\nabla \log\circ
\rho_i=\tilde{U}_i+W$, $1\leq i\leq k$. Then, for $i=1,\ldots,k$
we have that $E_i$ is umbilical with respect to $\<\;,\;\>$ with
mean curvature normal $H_i=U_i^{\perp_i}=\tilde{U}_i+W^{\perp_i}$.
Moreover, it follows from Lemma \ref{le:rhos} that $E_i^\perp$ is
also umbilical with mean curvature normal $\eta_i=W^i$. Therefore,
for any  $i\in \{1,\ldots,k\}$ we have
\begin{eqnarray}\label{eq:12} \<\nabla_{X_{\perp_i}}\eta_i,X_i\>
&=&\<\nabla_{X_{\perp_i}}W,X_i\>-\<\nabla_{X_{\perp_i}}W^{\perp_i},X_i\>\nonumber\\
&=&\<\nabla_{X_{\perp_i}}W,X_i\>-\<X_{\perp_i},W^{\perp_i}\>\<X_i,W^i\>.
\end{eqnarray}
On the other hand, using that $\tilde{U}_i\in E_0$ is a gradient vector field  we have
$$ \<\nabla_{X_i}\tilde{U}_i,X_{\perp_i}\>=\<\nabla_{X_{\perp_i}}\tilde{U}_i,X_i\>=\<X_{\perp_i},\tilde{U}_i\>\<W^i,X_i\>,$$
and hence
\begin{eqnarray} \label{eq:21}\<\nabla_{X_i}H_i,X_{\perp_i}\>&=&\<\nabla_{X_i}W,X_{\perp_i}\>-\<\nabla_{X_i}W^{i},X_{\perp_i}\>+\<\nabla_{X_i}\tilde{U}_i,X_{\perp_i}\>\nonumber
\\
                                                    &=&\<\nabla_{X_i}W,X_{\perp_i}\>-\<X_i,W^i\>\<H_i,X_{\perp_i}\>+\<X_{\perp_i},\tilde{U}_i\>\<W^i,X_i\>\nonumber\\
&=&\<\nabla_{X_i}W,X_{\perp_i}\>-\<X_i,W^i\>\<W^{\perp_i},X_{\perp_i}\>.
\end{eqnarray}
Since $W$ is a gradient vector field, we conclude from
(\ref{eq:12}) and (\ref{eq:21}) that (\ref{eq:hs}) holds, and thus
${\cal E}$ is a $CWP$-net with respect to $\<\;,\;\>$. Moreover,
if $\<\;,\;\>^\sim$ is a Riemannian product metric on $M$, then
also $E_0^\perp$ is umbilical with respect to $\<\;,\;\>$ with
mean curvature normal $\eta_0=W^{0}$. Therefore ${\cal E}$ is a
$CP$-net with respect to $\<\;,\;\>$.

   We now prove the converse. If ${\cal E}=(E_i)_{i=0,\ldots,k}$ is a $CWP$-net with respect to
   $\<\;,\;\>$ then it is also a $TP$-net by Lemma \ref{prop:cwp}, thus $\<\;,\;\>$ is a
   twisted product metric $\<\;,\;\>=\sum_{i=0}^k\rho_i^2\pi_i^*\<\;,\;\>_i$
   by Proposition \ref{prop:tpwp}.  Moreover, by  Lemma \ref{le:rhos} we have
   that $\rho_i/\rho_0$ only depends on $M_0$ and $M_i$ for $i=1,\ldots,
   k$.\vspace{1ex}\\
   CLAIM: There exist $\tilde{\va}_i\in C^\infty(M_i)$
   and $\tilde{\psi}_i\in C^\infty(M_0)$ such that $\rho_i/\rho_0=
   (\tilde{\va}_i\circ \pi_i)(\tilde{\psi}_i\circ \pi_0)$
   for  $i=1,\ldots, k$.\vspace{1ex}\\
   Assuming the claim, we conclude that $\<\;,\;\>$ is conformal to the warped product metric
$\<\;,\;\>^\sim=\pi_0^*\<\;,\;\>_0+\sum_{i=1}^k\psi_i^2\pi_i^*\<\;,\;\>^\sim_i$
with conformal factor $\rho_0$, where $\psi_i=\tilde{\psi}_i\circ
\pi_0$ and $\<\;,\;\>^\sim_i=\tilde{\va}^2_i\<\;,\;\>_i$ for
$i=1,\ldots,k$.\vspace{1ex}\\
{\it Proof of the claim:}  Let $U_i=-\nabla \log\circ\rho_i$.
   Then  $(U_i-U_0)=(U_i-U_0)^i+(U_i-U_0)^0$, and the claim is equivalent to $(U_i-U_0)^i$
   (and hence also $(U_i-U_0)^0$) being a gradient vector field.  Since $M$ is simply connected,
   this is in turn equivalent to $\<\nabla_{X}(U_i-U_0)^i,Y\>$ being
   symmetric in $X$ and $Y$. We now verify that this is indeed the case. First, using that $E_i^\perp$ is
   umbilical with mean curvature normal $\eta_i$ we have
   $$\<\nabla_{X_{\perp_i}}(U_i-U_0)^i,Y_{\perp_i}\>=
   -\<X_{\perp_i},Y_{\perp_i}\>\<\eta_i,(U_i-U_0)^i\>=
   \<\nabla_{Y_{\perp_i}}(U_i-U_0)^i,X_{\perp_i}\>.$$
 For $X=X_i$ and $Y=Y_i$, using that $E_i$ is umbilical with mean curvature normal $H_i$,
  the symmetry follows from
\begin{eqnarray*}\<\nabla_{X_i}(U_i-U_0)^i,Y_{i}\>&=&
\<\nabla_{X_i}(U_i-U_0),Y_{i}\>-\<\nabla_{X_i}(U_i-U_0)^0,Y_{i}\>\\
&=&\<\nabla_{X_i}(U_i-U_0),Y_{i}\>+\<X_i,Y_i\>\<H_i,(U_i-U_0)^0\>
\end{eqnarray*}
and the fact that $U_i-U_0$ is a gradient vector field. Finally,
we consider the case $X=X_{\perp_i}$ and $Y=X_i$. On one hand,
using that the mean curvature normal of $E_i^\perp$ is given by
$\eta_i=H_0^i=U_0^i$ by Lemma \ref{le:rhos}, and that
$U_i^{\perp_i}=H_i$, we have
\begin{eqnarray}\label{eq:sym1}\<\nabla_{X_{\perp_i}}(U_i-U_0)^i,X_i\>
&=&\<\nabla_{X_{\perp_i}}U_i,X_i\>-\<\nabla_{X_{\perp_i}}H_i,X_i\>-
\<\nabla_{X_{\perp_i}}\eta_i,X_i\>\nonumber\\
&=&\<\nabla_{X_{\perp_i}}U_i,X_i\>-\<X_{\perp_i},H_i\>\<\eta_i,X_i\>-
\<\nabla_{X_{\perp_i}}\eta_i,X_i\>.
\end{eqnarray}
On the other hand,
\begin{eqnarray}\label{eq:sym2}\<\nabla_{X_i}(U_i-U_0)^i,X_{\perp_i}\>&=&\<\nabla_{X_i}U_i,X_{\perp_i}\>
-\<\nabla_{X_i}H_i,X_{\perp_i}\>-\<\nabla_{X_i}\eta_i,X_{\perp_i}\>\nonumber\\
&=&\<\nabla_{X_i}U_i,X_{\perp_i}\>-\<\nabla_{X_i}H_i,X_{\perp_i}\>-\<X_i,\eta_i\>\<X_{\perp_i},H_i\>.
\end{eqnarray}
It follows from (\ref{eq:hs}) and the fact that $U_i$ is a
gradient vector field that the right-hand-sides of (\ref{eq:sym1})
and (\ref{eq:sym2}) coincide, and the proof of the claim is
completed.

    Now assume that ${\cal E}$ is a $CP$-net with respect to $\<\;,\;\>$, that is,
    also $E_0^\perp$ is umbilical.  We may assume that $k\geq 2$. It follows from Lemma \ref{le:rhos}
    that $\psi_i/\psi_1$ does not depend on $M_0$ for   $i=2,\ldots, k$, thus there exist
     $a_i\neq 0$ such that $\psi_i=a_i\psi_1$ for   $i=2,\ldots, k$.
     Therefore $\<\;,\;\>^\sim$ is conformal to the Riemannian product metric
     $\<\;,\;\>^*=\sum_{i=0}^k\pi_i^*\<\;,\;\>_i^*$ with conformal factor $\psi_1$,
     where $\<\;,\;\>^*_1=\<\;,\;\>^\sim_1$, $\<\;,\;\>^*_i=a_i^2\<\;,\;\>^\sim_i$
     for  $i=2,\ldots, k$ and
     $\<\;,\;\>^*_0=\tilde{\psi}_1^{-2}\<\;,\;\>_0$.\vspace{1ex}\qed

Arguing as in the proof of Theorem \ref{prop:qwp}, we obtain
from Proposition \ref{prop:conf}  the following conformal versions
of the local de Rham and Hiepko theorems.

\begin{theorem} \label{cor:drh} Let ${\cal E}=(E_i)_{i=0,\ldots,k}$ be a $CWP$-net
(resp., $CP$-net) on a Riemannian manifold $M$. Then for every
point $p\in M$ there exists a local product representation
$\psi\colon\;\Pi_{i=0}^k M_i\to U$ of ${\cal E}$ with $p\in
U\subset M$, which is a conformal diffeomorphism with respect to a
warped (resp., Riemannian) product metric on $\Pi_{i=0}^k M_i$.
\end{theorem}

 Given a connected product manifold $M=\Pi_{i=0}^k M_i$  endowed with a metric
   $\<\;,\;\>$ that is conformal to a warped product metric on
   $M$, we now investigate under what conditions  the
   subbundles
   $E_i$ and $E_i^\perp$ of the product net
   ${\cal E}=(E_i)_{i=0,\ldots,k}$ of $M$ are spherical with  respect to
   $\<\;,\;\>$ for a fixed $i\in \{1,\ldots, k\}$.

\begin{lemma} \label{le:sph} Let $M=\Pi_{i=0}^k M_i$ be a connected and simply connected
product manifold and let  ${\cal E}=(E_i)_{i=0,\ldots,k}$ be its
product net. Let
$\<\;,\;\>^\sim=\pi_0^*\<\;,\;\>^\sim_0+\sum_{i=1}^k
\tilde{\rho}_i^2\pi_i^*\<\;,\;\>^\sim_i$ be a warped product
metric on $M$ and  let $\<\;,\;\>=\va^2\<\;,\;\>^\sim$ be
conformal to $\<\;,\;\>^\sim$ with conformal factor $\va\in
C^\infty(M)$. Then for each fixed $i\in \{1,\ldots, k\}$ the
following assertions are equivalent:
\begin{itemize}
\item[{\em $(i)$}] Either $E_i$ or $E_i^\perp$ is (and hence both $E_i$
and $E_i^\perp$ are) spherical with respect to $\<\;,\;\>$;
\item[{\em $(ii)$}] The vector field $W=-\nabla \log\circ\va$ satisfies
\be\label{eq:ii}\<\nabla_{X_{\perp_i}}W,X_i\>=\<X_{\perp_i},W\>\<X_i,W\>;\ee
\item[{\em $(iii)$}] The conformal factor $\va\in
C^\infty(M)$ satisfies $\hess\,\va(X_i,X_{\perp_i})=0;$
\item[{\em $(iv)$}] $\tilde{\rho}^{-1}_i\tilde{W}^i$  is a gradient vector field, where $\tilde{W}=\nabla \va^{-1}$;
\item[{\em $(v)$}] There exist $\phi_i\in C^\infty(M_i)$ and
$\phi_{\perp_i}\in
C^\infty(M_{\perp_i}:=M_0\times\cdots\times\hat{M}_i\times\cdots\times
M_k)$ (the hat indicates that $M_i$ is missing) such that
$\va^{-1}=\tilde{\rho}_i(\phi_i\circ\pi_i)+\phi_{\perp_i}\circ
\pi_{\perp_i}$, where $\pi_{\perp_i}\colon\,M\to M_{\perp_i}$
denotes the canonical projection.
\end{itemize}
\end{lemma}
\proof  We have from the beginning of the proof of Proposition
\ref{prop:conf} that $E_i$ and $E_i^\perp$ are umbilical with
respect to $\<\;,\;\>$ with mean curvature normals
$H_i=\tilde{U}_i+W^{\perp_i}$ and $\eta_i=W^i$, respectively,
where $\tilde{U}_i =-\nabla \log\circ \tilde{\rho}_i\in E_0$. The
equivalence between $(i)$ and $(ii)$ then follows from
(\ref{eq:12}) and (\ref{eq:21}). Now, we have
\begin{eqnarray*}\nabla_XW&=&-\nabla_X(\va^{-1}\nabla\va)=
-X(\va^{-1})\nabla\va-\va^{-1}\nabla_X\nabla\va\\
&=&\va^{-2}\<\nabla\va,X\>\nabla\va-\va^{-1}\nabla_X\nabla\va\\
&=&\<W,X\>W-\va^{-1}\nabla_X\nabla\va. \end{eqnarray*}
Thus
$$\<\nabla_{X_{\perp_i}}W,X_i\>-\<X_{\perp_i},W\>\<X_i,W\>=-\va^{-1}\hess\,\va(X_{\perp_i},X_i),$$
and the equivalence between $(ii)$ and $(iii)$ follows. We now
prove that $(ii)$ is equivalent to the symmetry of
$\<\nabla_X\tilde{\rho}_i^{-1}\tilde{W}^i,Y\>$ with respect to $X$
and $Y$. Since $M$ is simply connected, this implies the
equivalence between $(ii)$ and $(iv)$. From $\tilde{W}=\va^{-1}W$
we obtain
$$\nabla_{X}\tilde{W}=\va^{-1}(\<W,X\>W+\nabla_X W),$$
and hence
\begin{eqnarray}\label{eq:calc}\<\nabla_{X_{\perp_i}} \tilde{W}^i, X_i\>
&=&\<\nabla_{X_{\perp_i}} \tilde{W}, X_i\>-\va^{-1}\<W^{\perp_i},X_{\perp_i}\>\<W^i,X_i\>\nonumber\\
&=&\va^{-1}\<\nabla_{X_{\perp_i}} W, X_i\>,
\end{eqnarray}
where in the first equality we have used that $E_i^\perp$ is
umbilical with mean curvature normal $W^i$. Since
$$X_{\perp_i}(\tilde{\rho}_i^{-1})=\tilde{\rho}_i^{-1}X_{\perp_i}(-\log\circ
\tilde{\rho}_i)=\tilde{\rho}_i^{-1}\<\tilde{U}_i,X_{\perp_i}\>,$$
we obtain from (\ref{eq:calc}) that
\begin{eqnarray}\label{eq:g1}\<\nabla_{X_{\perp_i}}
\tilde{\rho}_i^{-1}\tilde{W}^i, X_i\>
&=&\tilde{\rho}_i^{-1}\<\tilde{U}_i,X_{\perp_i}\>\<\tilde{W}^i,X_i\>+
\tilde{\rho}_i^{-1}\va^{-1}\<\nabla_{X_{\perp_i}}W,X_i\>\nonumber\\
&=&\tilde{\rho}_i^{-1}\va^{-1}(\<\tilde{U}_i,X_{\perp_i}\>\<W^i,X_i\>+
\<\nabla_{X_{\perp_i}}W,X_i\>).
\end{eqnarray}
 On the other hand, since $E_i$ is umbilical with mean curvature
 normal $H_i=\tilde{U}_i+W^{\perp_i}$ we have
\begin{eqnarray*}
\<\nabla_{X_i}
\tilde{W}^i,X_{\perp_i}\>&=&\<X_i,\tilde{W}^i\>\<X_{\perp_i},H_i\>\\
&=&\va^{-1}(\<X_i,W^i\>\<X_{\perp_i},W^{\perp_i}\>+
\<X_i,W^i\>\<X_{\perp_i},\tilde{U}_i\>). \end{eqnarray*} Using
that $\tilde{\rho}_i$ only depends on $M_0$, we obtain
\be\label{eq:g2}
\<\nabla_{X_i}\tilde{\rho}_i^{-1}\tilde{W}^i,X_{\perp_i}\>
=\tilde{\rho}_i^{-1}\va^{-1}(\<X_i,W^i\>\<X_{\perp_i},W^{\perp_i}\>+
\<X_i,W^i\>\<X_{\perp_i},\tilde{U}_i\>).
 \ee
Comparing (\ref{eq:g1}) and (\ref{eq:g2}) implies that
$\<\nabla_{X_{\perp_i}} \tilde{\rho}_i^{-1}\tilde{W}^i,
X_i\>=\<\nabla_{X_i} \tilde{\rho}_i^{-1}\tilde{W}^i,X_{\perp_i}\>$
 is equivalent to $(ii)$.
 Since
$$\<\nabla_{X_{\perp_i}} \tilde{\rho}_i^{-1}\tilde{W}^i,Y_{\perp_i}\>=
-\tilde{\rho}_i^{-1}\<X_{\perp_i},Y_{\perp_i}\>\<W^i,\tilde{W}^i\>=\<\nabla_{Y_{\perp_i}}
\tilde{\rho}_i^{-1}\tilde{W}^i,X_{\perp_i}\>,$$ where we have used
again that $E_i^\perp$ is umbilical with mean curvature normal
$W^i$, and
\begin{eqnarray*} \<\nabla_{X_i} \tilde{\rho}_i^{-1}\tilde{W}^i,Y_i\>
&=&\tilde{\rho}_i^{-1}(\<\nabla_{X_i} \tilde{W},Y_i\>-\<\nabla_{X_i} \tilde{W}^{\perp_i},Y_i\>)\\
&=&\tilde{\rho}_i^{-1}(\<\nabla_{X_i}
\tilde{W},Y_i\>+\<X_i,Y_i\>\<H_i,\tilde{W}^{\perp_i}\>)=\<\nabla_{Y_i}
\tilde{\rho}_i^{-1}\tilde{W}^i,X_i\>,
\end{eqnarray*}
because $\tilde{W}$ is a gradient vector field and $E_i$ is
umbilical with mean curvature normal $H_i$, the proof of the
equivalence between $(ii)$ and $(iv)$ is completed.

    Finally, if $(iv)$ holds  let $\phi_i\in C^\infty(M_i)$ be such that
$\tilde{\rho}_i^{-1}\tilde{W}^i=\nabla(\phi_i\circ\pi_i)$. Then
$\nabla(\va^{-1}-\tilde{\rho}_i(\phi_i\circ\pi_i))\in
E_{\perp_i}$, hence there exists $\phi_{\perp_i}\in
C^\infty(M_{\perp_i})$  such that
$\va^{-1}-\tilde{\rho}_i(\phi_i\circ\pi_i)=\phi_{\perp_i}\circ
\pi_{\perp_i}$. Conversely, assuming $(v)$ we obtain
$\tilde{W}^i=\tilde{\rho}_i\nabla(\phi_i\circ\pi_i)$, thus $(iv)$
holds. \qed

\begin{remark} \label{re:scon}{\em The assumption that $M$ is
simply connected in Lemma \ref{le:sph} may be dropped if the
assertions in $(iv)$ and $(v)$ are required to hold only locally.}
\end{remark}

\begin{proposition} \label{cor:sph} Let $M=\Pi_{i=0}^k M_i$ be a connected and simply connected
product manifold and let  ${\cal E}=(E_i)_{i=0,\ldots,k}$ be its
product net.  Then  the following assertions on a Riemannian
metric $\<\;,\;\>$ on $M$ are equivalent:
\begin{itemize}
\item[{\em $(i)$}] $E_i$ and $E_i^\perp$ are spherical for $i=1,\dots,k$
 (resp., $E_i$ and $E_i^\perp$ are spherical for $i=1,\dots,k$ and, in addition, $E_0^\perp$
 is spherical) with respect to $\<\;,\;\>$;
\item[{\em $(ii)$}] $\<\;,\;\>=\va^2\<\;,\;\>^\sim$ is conformal to a warped product metric
$$\<\;,\;\>^\sim=\pi_0^*\<\;,\;\>^\sim_0+\sum_{i=1}^k
\tilde{\rho}_i^2\pi_i^*\<\;,\;\>^\sim_i$$ (resp., a Riemannian
product metric $\<\;,\;\>^\sim$) on $M$ with conformal factor
$\va\in C^\infty(M)$ given by $\va^{-1}=\phi_0\circ
\pi_0+\sum_{i=1}^k \tilde{\rho}_i(\phi_i\circ \pi_i)$ (resp.,
$\va^{-1}=\sum_{i=0}^k \phi_i\circ \pi_i$), where $\phi_i\in
C^\infty(M_i)$ for $i=0,\ldots,k$. \end{itemize}
\end{proposition}
\proof Assuming $(i)$, we have from Proposition \ref{prop:conf}
that $\<\;,\;\>=\va^2\<\;,\;\>^\sim$ is conformal to a warped
product metric
$\<\;,\;\>^\sim=\pi_0^*\<\;,\;\>^\sim_0+\sum_{i=1}^k
\tilde{\rho}_i^2\pi_i^*\<\;,\;\>^\sim_i$ (resp., a Riemannian
product metric $\<\;,\;\>^\sim$) on $M$ with conformal factor
$\va\in C^\infty(M)$. Moreover, by Lemma \ref{le:sph} there exist
$\phi_i\in C^\infty(M_i)$ such that
$\tilde{\rho}_i^{-1}(\nabla\va^{-1})^i=\nabla(\phi_i\circ \pi_i)$
(resp., $(\nabla\va^{-1})^i=\nabla(\phi_i\circ \pi_i)$) for
$i=1,\ldots,k$. Therefore
$\nabla(\va^{-1}-\sum_{i=1}^k\tilde{\rho}_i(\phi_i\circ \pi_i))\in
E_0$ (resp., $\nabla(\va^{-1}-\sum_{i=1}^k\phi_i\circ \pi_i)\in
E_0$), which implies $(ii)$. Conversely, if $(ii)$ holds for a
warped product metric $\<\;,\;\>^\sim$ on $M$, then $(i)$ holds by
the fact that $(v)$ implies $(i)$ in Lemma \ref{le:sph} for
$i=1,\ldots,k$. Moreover, if $\<\;,\;\>^\sim$ is a Riemannian
product metric, then $\E$ is a $CP$-net by Proposition
\ref{prop:conf}, thus also $E_0^\perp$ is umbilical. But
$E_i^\perp$ being spherical for $i=1,\ldots, k$ implies that $E_0$
is spherical by Lemma \ref{prop:cwp}, and hence also $E_0^\perp$
is spherical by Proposition \ref{prop:sph}-$(iii)$.
\vspace{1ex}\qed

Arguing again as in the proof of Theorem \ref{prop:qwp}, we
obtain from Proposition \ref{cor:sph}  the following local
decomposition theorem of de Rham-type.

\begin{theorem} \label{cor:sph2} Let ${\cal E}=(E_i)_{i=0,\ldots,k}$ be an orthogonal net
 on a Riemannian manifold $M$ such that $E_i$ and $E_i^\perp$ are
spherical for  $i=1,\ldots, k$ (resp., $E_i$ and $E_i^\perp$ are
spherical for  $i=1,\ldots, k$ and, in addition, $E_0^\perp$ is
spherical). Then for every point $p\in M$ there exists a local
product representation $\psi\colon\;\Pi_{i=0}^k M_i\to U$ of
${\cal E}$ with $p\in U\subset M$, which is a conformal
diffeomorphism with respect to a warped (resp., Riemannian)
product metric on $\Pi_{i=0}^k M_i$ such that the conformal factor
$\va\in C^{\infty}(M)$ of $\psi$ is given by $\va^{-1}=\phi_0\circ
\pi_0+\sum_{i=0}^k \tilde{\rho}_i(\phi_i\circ \pi_i)$ (resp.,
$\va^{-1}=\sum_{i=0}^k \phi_i\circ \pi_i$), where $\phi_i\in
C^\infty(M_i)$ for $i=0,\ldots,k$.
\end{theorem}

 \begin{remark} \label{re:classical}{\em  Classically, a one-parameter family of curves in a surface $S\subset \R^3$ is said to be an
  {\it isothermal family\/} if it is the family of $u$-coordinate curves of an
  isothermic coordinate system $(u,v)$ on $S$. In this setting, Theorem \ref{cor:drh}
  reduces to a well-known characterization of isothermal
  families of curves in terms of their geodesic curvatures and the
  geodesic curvatures of their orthogonal trajectories
  (cf. \cite{da}, vol. III, p. 154, eq.(36)). Theorem \ref{cor:drh} and
  Proposition \ref{prop:sph} also generalize the facts that orthogonal families
  of curves with constant geodesic curvature are isothermal, and that the curves of an
  isothermal family must have constant geodesic curvature if the same holds for
  their orthogonal trajectories (cf. \cite{da}, vol. III, p. 154). Finally,
  Theorem \ref{cor:sph2} extends the characterization  of the
  first fundamental form of a surface that admits two orthogonal
  families of curves with constant geodesic
  curvature (cf. \cite{bi}, p. 368, eq.~(38)).}
  \end{remark}\vspace{2ex}

\noindent {\bf\large  \S 5. Codazzi tensors. }\vspace{2ex}

     In this section we apply our results to study Riemannian manifolds
     that carry a Codazzi tensor with exactly two distinct eigenvalues everywhere.
     Recall that a symmetric tensor $\Phi$ is said to be a {\it Codazzi tensor\/}  if
$(\nabla_X\Phi)Y=(\nabla_Y\Phi) X$ for all $X,Y\in TM$, where
$(\nabla_X\Phi)Y=\nabla_X\Phi Y-\Phi\nabla_X Y$.  We start with
the following basic result due to Reckziegel \cite{re}, a short
proof of which is included for the sake of completeness.

\begin{proposition} \label{prop:codtensor} Let $\Phi$ be a Codazzi tensor on a
Riemmanian manifold $M$, and let $\lambda\in C^\infty(M)$ be an eigenvalue of
$\Phi$ such that $E_\lambda=\ker(\lambda I-\Phi)$ has constant rank $r$.
Then the following hold:
\begin{itemize}
\item[{\em (i)}] $E_\lambda$ is an umbilical distribution with mean curvature normal
$\eta$ given by \be\label{eq:mcn}(\lambda I-\Phi)\,\eta=(\nabla
\,\lambda)_{E_\lambda^\perp}.\ee
\item[{\em (ii)}] If $r\geq 2$ then $\lambda$ is
constant along $E_\lambda$.
\item[{\em (iii)}] If $\lambda$ is constant along $E_\lambda$ then $E_\lambda$ is spherical.
\end{itemize}
\end{proposition}
\proof Let $T\in E_\lambda$ and $X\in TM$. Taking the inner
product of both sides of $(\nabla_X\Phi)T = (\nabla_T\Phi)X$ with
$S\in E_\lambda$ yields \be\label{eq:lamb}(\lambda I-\Phi)\nabla_T
S=\<T,S\>\nabla\,\lambda -T(\lambda)S.\ee Since $\lambda I-\Phi$
vanishes on $E_\lambda$, we have that $(\lambda I-\Phi)\nabla_T
S=(\lambda I-\Phi)(\nabla_T S)_{E_\lambda^\perp}\in
E_\lambda^\perp$. Therefore, comparing the components in
 $E_\lambda^\perp$ of both sides of (\ref{eq:lamb})
yields \be\label{eq:lamb2}(\lambda I-\Phi)(\nabla_T
S)_{E_\lambda^\perp}=\<T,S\>(\nabla
\,\lambda)_{E_\lambda^\perp},\ee and $(i)$ follows. Now assume
that $r\geq 2$. Since the left-hand-side of (\ref{eq:lamb}) is in
$E_\lambda^\perp$, it follows that $\<T,S\>(\nabla
\,\lambda)_{E_\lambda}=T(\lambda)S.$ Then, for any $T\in
E_\lambda$, choosing  $0\neq S\in E_\lambda$ orthogonal to $T$
yields $T(\lambda)=0$. To prove $(iii)$, from $T(\lambda)=0$ for
all $T\in E_\lambda$  we have
\begin{eqnarray*}  \<\nabla_T \eta, (\lambda I-\Phi)X\>&=& T\<(\lambda I-\Phi)\eta,X\>-\<\eta,\nabla_T  (\lambda I-\Phi)X\>\\
&=& TX(\lambda)-\lambda\<\nabla_T X,\eta\>+\<\nabla_T \Phi
X,\eta\>.
\end{eqnarray*}
Using that $\nabla_T \Phi X =\nabla_X \Phi T  -\Phi \nabla_X T +
\Phi \nabla_T X$ we obtain
$$
\<\nabla_T \Phi X,\eta\>=\<(\lambda I-\Phi)\eta,\nabla_X
T\>+\<\Phi \eta,\nabla_T X\>.
$$
Therefore
$$\<\nabla_T \eta, (\lambda I-\Phi)X\>=TX(\lambda)-\<(\lambda I-\Phi)\eta,[T,X]\>=TX(\lambda)-[T,X](\lambda)=0.\qed$$

Here and in the sequel, writing a vector subbundle as a subscript
of a vector field indicates taking its component in that
subbundle. From now on we consider Codazzi tensors  with exactly
two distinct eigenvalues
 $\lambda$ and $\mu$ everywhere, and always denote sections of the corresponding eigenbundles $E_\lambda$ and $E_\mu$ by $X$ and $Y$, respectively.

\begin{theorem} \label{thm:codtensor} Let $M$ be a connected Riemannian manifold and let $\Phi$ be a
Codazzi tensor on $M$ with exactly two distinct eigenvalues
$\lambda$ and $\mu$ everywhere. Let $E_\lambda$ and $E_\mu$ be the
corresponding eigenbundles. Then, for every point $p\in M$ there
exists a local product representation $\psi\colon\;M_0\times
M_1\to U$ of $(E_\lambda, E_\mu)$ with $p\in U\subset M$, which is
an isometry with respect to a twisted product metric $\<\;,\;\>$
on $M_0\times M_1$. Moreover,
\begin{itemize}
\item[{\em (i)}]  $\<\;,\;\>$ is conformal to a Riemannian product metric  if and only if
\be\label{eq:cpnet} 2\beta X(\alpha)Y(\alpha)+\alpha
X(\alpha)Y(\beta)+\alpha Y(\alpha)X(\beta)-\alpha\beta\hess\,
\alpha(X,Y)=0, \ee where
$\displaystyle{\alpha=\frac{1}{2}(\mu+\lambda)}$ and
$\beta=\displaystyle{\frac{\mu-\lambda}{\mu+\lambda}}$.
\item[{\em (ii)}] If either
\be\label{eq:mulambda1} 2X(\mu)Y(\mu)-X(\mu)Y(\lambda)+(\lambda-\mu)\hess \mu(X,Y)=0\ee
or
\be\label{eq:mulambda2} 2X(\lambda)Y(\lambda)-X(\mu)Y(\lambda)-(\lambda-\mu)\hess \lambda(X,Y)=0,\ee
 then (shrinking $U$ if necessary) the assertions below are equivalent:
\begin{itemize}
\item[{\em (a)}]  $\<\;,\;\>$ is conformal to a Riemannian product metric on $M_0\times M_1$;
\item[{\em (b)}]  $\<\;,\;\>=\va^2\<\;,\;\>^\sim$, where $\<\;,\;\>^\sim$ is a Riemannian product metric
 on $M_0\times M_1$ and $\va^{-1}=\phi_0\circ
\pi_0+\phi_1\circ \pi_1$ for $\phi_i\in C^{\infty}(M_i)$, $0\leq
i\leq 1$.
\item[{\em (c)}] Both (\ref{eq:mulambda1}) and (\ref{eq:mulambda2}) hold.
\end{itemize}
\item[{\em (iii)}] Equation (\ref{eq:mulambda1}) (resp., (\ref{eq:mulambda2}))
is satisfied if $\lambda$ (resp., $\mu$) is constant along
$E_\lambda$ (resp., $E_\mu$), in particular if $E_\lambda$ (resp.,
$E_\mu$) has rank at least two. Therefore, the assumption in
$(ii)$ is always satisfied if $n\geq 3$. Moreover, if both
$\lambda$ and $\mu$ are constant along the corresponding
eigenbundles, then the functions $\phi_i\in C^{\infty}(M_i)$,
$0\leq i\leq 1$, in $(ii)\!-\!(b)$ can be chosen so that
$\psi_*^{-1}\Phi\psi_*=(\phi_1\circ \pi_1)\Pi_0-(\phi_0\circ
\pi_0)\Pi_1,$ where $\Pi_0$ and $\Pi_1$ denote the orthogonal
projections onto the subbundles $E_0$ and $E_1$, respectively,  of
the product net $(E_0, E_1)$ of $M_0\times M_1$.
\end{itemize}
\end{theorem}
\proof  By Proposition \ref{prop:codtensor},  $E_\lambda$ and
$E_\mu$ are umbilical  with mean curvature normals
\be\label{eq:pns} \eta=(\lambda-\mu)^{-1}(\nabla
\,\lambda)_{E_\mu}\;\;\;\;\mbox{and}\;\;\;\;\zeta=(\mu
-\lambda)^{-1}(\nabla\, \mu)_{E_\lambda}. \ee Thus $(E_\lambda,
E_\mu)$ is a $TP$-net and the first assertion follows from
Theorem~\ref{cor:isodrh}. Now we have \be
\label{eq:s1a}\<\nabla_{X}\eta,Y\>=X(\lambda-\mu)^{-1}Y(\lambda)+
(\lambda-\mu)^{-1}\<\nabla_{X}(\nabla \lambda)_{E_\mu},Y\>.\ee
Using that $E_\lambda$ is umbilical with mean curvature normal
$\eta$ we obtain \begin{eqnarray} \label{eq:s1b}
\<\nabla_{X}(\nabla \lambda)_{E_\mu},Y\>&=&\<\nabla_{X}\nabla
\lambda,Y\>-\<\nabla_{X}(\nabla
\lambda)_{E_\lambda},Y\>\nonumber\\
&=&
\hess\lambda(X,Y)-\<X,\nabla\lambda\>\<\eta,Y\>\nonumber\\
&=&\hess\lambda(X,Y)-(\lambda-\mu)^{-1}X(\lambda)Y(\lambda).
\end{eqnarray}
Substituting (\ref{eq:s1b}) into (\ref{eq:s1a}) yields
\be\label{eq:s1}  \<\nabla_{X}\eta,Y\>=
 -(\lambda-\mu)^{-2}(2X(\lambda)Y(\lambda)-X(\mu)Y(\lambda)-(\lambda-\mu)\hess
 \lambda(X,Y)).
 \ee
A similar computation using that $E_\mu$ is umbilical with mean
curvature normal $\zeta$ gives
 \be
\label{eq:s2}\<\nabla_{Y}\zeta,X\>=
-(\lambda-\mu)^{-2}(2X(\mu)Y(\mu)-X(\mu)Y(\lambda)+(\lambda-\mu)\hess
\mu(X,Y)). \ee Hence \be\label{eq:nzeta}
\begin{array}{l}(\lambda-\mu)^2(\<\nabla_{X}\eta,Y\>-\<\nabla_{Y}\zeta,X\>)=\vspace{1.5ex}\\
\hspace*{6ex}-2X(\lambda)Y(\lambda)+2X(\mu)Y(\mu)+(\lambda-\mu)\hess (\lambda+\mu)(X,Y).
\end{array}
\ee Now, a straightforward computation using that
$\lambda=\alpha(1-\beta)$ and $\mu=\alpha(1+\beta)$ yields
\be\label{eq:aux} X(\lambda)Y(\lambda)-X(\mu)Y(\mu)=-4\beta
X(\alpha)Y(\alpha)-2\alpha X(\alpha)Y(\beta)-2\alpha
Y(\alpha)X(\beta). \ee
Substituting (\ref{eq:aux}) into (\ref{eq:nzeta}) and observing
that $\lambda+\mu=2\alpha$ and $\lambda-\mu=-2\alpha\beta$, we
obtain that $(E_\lambda, E_\mu)$ is a $CP$-net if and only if
(\ref{eq:cpnet}) holds. The assertion in $(i)$ now follows from
Proposition~\ref{prop:conf}.
We now prove $(ii)$. By (\ref{eq:s1}) and (\ref{eq:s2}), we have that  (\ref{eq:mulambda1}) and (\ref{eq:mulambda2}) are equivalent to $E_\mu$ and $E_\lambda$ being spherical, respectively. Therefore, if either (\ref{eq:mulambda1}) or (\ref{eq:mulambda2}) holds, then the following are equivalent by  Proposition~\ref{prop:sph}:\vspace{1ex}\\
$(a')$ $(E_\lambda, E_\mu)$ is a $CP$-net;$\;\;\;(b')$  $E_\lambda$ and $E_\mu$ are spherical;$\;\; \;$$(c)$ Both (\ref{eq:mulambda1}) and (\ref{eq:mulambda2}) hold.\vspace{1ex}\\
The proof of $(ii)$ is completed by observing that  $(a')$ is
equivalent to $(a)$ by Proposition~\ref{prop:conf} and that $(b')$
is equivalent to  $(b)$ by Proposition \ref{cor:sph} (shrinking $U$
if necessary).

 If  $\lambda$ (resp., $\mu$) is constant along $E_\lambda$ (resp., $E_\mu$),
 which is always the case by  Proposition \ref{prop:codtensor}-$(ii)$
 if $E_\lambda$ (resp., $E_\mu$) has rank at least two, then $E_\lambda$
 (resp., $E_\mu$) is spherical by  Proposition \ref{prop:codtensor}-$(iii)$. Therefore the first assertion in
 $(iii)$ follows from (\ref{eq:s1}) (resp., (\ref{eq:s2})). Now assume
 that both $\lambda$ and $\mu$ are constant along the
 corresponding eigenbundles. Denote  $\tilde{\phi}_1=\lambda\circ \psi$ and
$\tilde{\phi}_0=\mu\circ \psi$ the eigenvalues of
$\psi_*^{-1}\Phi\psi_*$. Then the eigenbundles of $\tilde{\phi}_0$
and $\tilde{\phi}_1$ are $E_0$ and $E_1$, respectively, thus we
have from (\ref{eq:pns}) that the mean curvature normals of $E_0$
and $E_1$ are given, respectively, by
$$\eta=(\tilde{\phi}_1-\tilde{\phi}_0)^{-1}\nabla \,\tilde{\phi}_1\;\;\;\;\mbox{and}\;\;\;\;\zeta=(\tilde{\phi}_0-\tilde{\phi}_1)^{-1}\nabla\, \tilde{\phi}_0.$$
On the other hand, since $(E_0, E_1)$  is the product net of
$M_0\times M_1$, we have
$$\eta=-(\nabla\, \log\circ \va)_{E_1}\;\;\;\;\mbox{and}
\;\;\;\;\zeta=-(\nabla\, \log\circ \va)_{E_0}.$$ It follows that
$\nabla \, \log\circ
(\tilde{\phi}_1-\tilde{\phi}_0)=-\nabla\,\log\circ \va$, thus
$\va^{-1}=A(\tilde{\phi}_1-\tilde{\phi}_0)$ for some $A\neq 0$. By
rescaling the metric $\<\;,\;\>$ we may assume that $A=1$. The
proof is completed by defining $\phi_i\in C^{\infty}(M_i)$, $0\leq
i\leq 1$, such that $\tilde{\phi}_0=-\phi_0\circ \pi_0$ and
$\tilde{\phi}_1=\phi_1\circ \pi_1$.\qed



 \begin{remark} \label{re:isot}{\em In case  $\Phi$ is the shape operator of a
 surface in $\R^3$ then (\ref{eq:cpnet})
reduces to a criterion for the surface to be isothermic in terms
of its principal curvatures. In particular, since it is clearly
satisfied if $\lambda+\mu=2\alpha$ is constant on $M$, it implies
the well-known fact that surfaces with constant mean curvature are
isothermic surfaces.} \end{remark}

\begin{corollary}\label{prop:codazzi23} Let $M$ be a connected Riemannian manifold and let $\Phi$ be a Codazzi tensor on $M$ with exactly two distinct eigenvalues $\lambda$ and $\mu$ everywhere. Let $E_\lambda$ and $E_\mu$ be the corresponding eigenbundles. Assume that $\mu$ is constant along $E_\mu$  and that $\lambda=h(\mu)$ for some smooth real function $h$. Then one of the following possibilities holds:
\begin{itemize}
\item[{\em (i)}] $\lambda$ is constant along $E_\lambda$ and for every point $p\in M$ there exists a local product representation $\psi\colon\;M_0\times M_1\to U$ of $(E_\lambda, E_\mu)$ with $p\in U\subset M$, which is an isometry with respect to a Riemannian product metric $\<\;,\;\>$ on $M_0\times M_1$. Moreover, $\psi_*^{-1}\Phi\psi_*=A_0\Pi_0+A_1\Pi_1$ for some $A_0,A_1\in \R$.
\item[{\em (ii)}] for every point $p\in M$ there exists a local product representation $\psi\colon\;I\times M_1\to U$ of $(E_\lambda, E_\mu)$ with $p\in U\subset M$, where $I\subset \R$ is an open interval, which is an isometry with respect to a warped  product metric $\<\;,\;\>$ on $I\times M_1$ with warping function $\sigma=\tilde{\sigma}\circ \pi_0$, where $\tilde{\sigma}\in C^\infty(I)$. Moreover, $\psi_*^{-1}\Phi\psi_*=h(\tilde{\mu}\circ\pi_0)\Pi_0+(\tilde{\mu}\circ\pi_0)\Pi_1$, where $\tilde{\mu}\in C^\infty(I)$ is determined (up to a constant) by
\be\label{eq:tilmu}\int \frac{d\tilde{\mu}}{h(\tilde{\mu})- \tilde{\mu}}= \log\tilde{\sigma}.\ee
\end{itemize}
\end{corollary}
\proof Since $\mu$ is constant along $E_\mu$ then $(\nabla \,\mu)_{E_\mu}=0$, and hence also $(\nabla \,\lambda)_{E_\mu}=0$ by the assumption that $\lambda=h(\mu)$ for some smooth real function $h$. It follows from  (\ref{eq:pns}) that $E_\lambda$ is totally geodesic. If also $\lambda$ is constant along $E_\lambda$, we obtain in a similar way that $E_\mu$ is totally geodesic,  thus
$(i)$ follows from  the local de Rham theorem and the fact that
both $\lambda$ and $\mu$ are now constant on $M$. Otherwise, by
Proposition \ref{prop:codtensor}-$(ii)$ the distribution
$E_\lambda$ must have rank one. Since $E_\mu$ is spherical by
Proposition \ref{prop:codtensor}-$(iii)$, then $(E_\lambda,
E_\mu)$ is a $WP$-net. Thus, the first assertion in $(ii)$ holds
by Theorem~\ref{cor:isodrh}. Finally, let $E_0$ and $E_1$ be the
eigenbundles of $\psi_*^{-1}\Phi\psi_*$ correspondent to its
eigenvalues $\lambda\circ \psi$ and $\mu\circ \psi$, respectively.
Since $(E_0,E_1)$ is the product net of  $I\times M_1$, the mean
curvature normal of $E_1$ is given, on one hand, by
$\zeta=-\nabla\,\log\circ \sigma$, and on the other hand by
$\zeta=(\mu\circ \psi-\lambda\circ \psi)^{-1}(\nabla\,\mu\circ
\psi)_{E_0}=(\mu\circ \psi-h(\mu\circ \psi))^{-1}\nabla\,\mu\circ
\psi$, because $(\nabla\,\mu\circ \psi)_{E_1}=0$. Writing
$\mu\circ \psi=\tilde{\mu}\circ \pi_0$ for $\tilde{\mu}\in
C^\infty(I)$, we conclude that $\tilde{\mu}$ and $\tilde{\sigma}$
are related by (\ref{eq:tilmu}).\qed

\begin{remarks} \label{re:derd}{\em $(i)$ Codazzi tensors with  two eigenvalues have
also been studied in \cite{de} and \cite{ds}. In \cite{ds} it is
shown that the existence of such a tensor on a Riemannian manifold
imposes some restrictions on its curvature tensor. One can check
that such restrictions follow from the form of the curvature
tensor of a twisted
product (cf. formula ($3$) in \cite{mrs}).  Corollary \ref{prop:codazzi23} is proved in \cite{de} for Codazzi tensors with constant trace.\vspace{1ex}\\
$(ii)$ Part $(ii)$ of Corollary \ref{prop:codazzi23}
can be regarded as an intrinsic version of Theorem $4.2$ in \cite{dcd},
which states that a  hypersurface of dimension $n\geq 3$ of a space form with two
distinct principal curvatures $\lambda$ and $\mu$, one of which has multiplicity one,
must be a rotation hypersurface whenever $\lambda=h(\mu)$ for some smooth real function $h$.
In fact, that result follows from Corollary \ref{prop:codazzi23} together with N\"olker's
decomposition theorem for isometric immersions of warped products
into space forms \cite{nol}.}\end{remarks}

{\renewcommand{\baselinestretch}{1} \hspace*{-25ex}\begin{tabbing}
\indent  \=    Universidade Federal de S\~ao Carlos \\
\indent  \= Via Washington Luiz km 235 \\
\> 13565-905 -- S\~ao Carlos -- Brazil \\
\> email: tojeiro@dm.ufscar.br \\
\end{tabbing}}

\begin{thebibliography}{lbllll}

\bibitem[{\bf $1$}]{bi}  BIANCHI, L.:  {\it Lezioni di Geometria Differenziale\/}, Bologna, 1927.

\bibitem[{\bf $2$}]{bu} BURSTALL, F.: Isothermic surfaces: conformal geometry,
Clifford algebras and integrable systems. Preprint,
math-DG/0003096.


\bibitem[{\bf $3$}]{da} DARBOUX, G.: {\it Le\c cons sur la th\'eorie des surfaces\/} (Reprinted by Chelsea Pub. Co., 1972), Paris 1914.

\bibitem[{\bf $4$}]{Rh} DE RHAM, G.: Sur la r\'educibilit\'e d'un espace de Riemann.  {\it Comm. Math. Helv.\/} ~{\bf 26} (1952), 328--344.

\bibitem[{\bf $5$}]{de} DERDZI\'NSKI, A.: Some remarks on the local structure of Codazzi tensors. Global differential geometry and global analysis (Berlin, 1979), pp.251--255. Lecture Notes in Math. {\bf 838}, Springer, Berlin-New York 1981.


\bibitem[{\bf $6$}]{ds} DERDZI\'NSKI, A.; SHEN, C.: Codazzi tensor fields, curvature and Pontryagin forms. {\it Proc. London Math. Soc.\/} ~{\bf 47} (1983), 15--26.


\bibitem[{\bf $7$}]{dcd} do CARMO, M.; DAJCZER, M.: Rotation hypersurfaces in spaces of constant curvature. {\it Trans. Amer. Mah. Soc.} {\bf 277} (1983), 685-709.


\bibitem[{\bf $8$}]{hi} HIEPKO, S.: Eine innere Kennzeichnung der verzerrten Produkte.  {\it Math. Ann.\/} ~{\bf 241} (1979), 209--215.



\bibitem[{\bf $9$}]{mrs} MEUMERTZHEIM, M.; RECKZIEGEL, H.; SCHAAF. M.: Decomposition of twisted and warped product nets.
{\it Result. Math.\/} {\bf 36} (1999), 297--312.

\bibitem[{\bf $10$}]{nol} N\"OLKER, S.:  Isometric immersions of warped products. {\it Diff. Geom. Appl.\/}~{\bf 6} (1996), 31--50.

\bibitem[{\bf $11$}]{rs} RECKZIEGEL, H.; SCHAAF. M.: De Rham decomposition of netted manifolds
{\it Result. Math.\/} {\bf 35} (1999), 175--191.

\bibitem[{\bf $12$}]{re} RECKZIEGEL, H.: Krummungsflachen von isometrischen immersionen in raume konstanter krummung.
{\it Math. Ann.\/} {\bf 223} (1976), 169--181.

\bibitem[{\bf $13$}]{to} TOJEIRO, R.: Isothermic submanifolds of Euclidean space.
Preprint.
\end{thebibliography}
\end{document}